\documentclass[11pt, a4]{article}

\usepackage{amsmath,amsfonts,amssymb,latexsym}

\usepackage{enumerate} 

\usepackage[pdftex]{graphicx,color}

\newenvironment{pf}{\par\textsc{Proof.}}{\medskip}
\newtheorem{thm}{Theorem}
\newtheorem{prop}[thm]{Proposition}
\newtheorem{lem}[thm]{Lemma}

\def\qed{\relax\ifmmode\hskip2em \Box 
\else\unskip\nobreak\hskip1em $\Box$\fi}

\title{Periodicity of power Fibonacci sequences\\
       modulus a Fibonacci number}
      
\author{Josep M. Brunat and Joan-C. Lario}

\date{April 1, 2022}

\begin{document}

\maketitle

\begin{abstract}

\noindent Let ${\mathcal F}=(F_i:i\ge 0)$ be the sequence of Fibonacci numbers, and $j$ and $e$ be non negative integers. We study the periodicity of the power Fibonacci sequences ${\mathcal F}^e(F_j)=(F_i^e\pmod{F_j}: i\ge 0)$. It is shown that for every $j,e\ge 1$ the sequence ${\mathcal F}^e(F_j)$ is periodic and its periodicity is computed. The result was previously known for
${\mathcal F}(F_j)$; that is, for $e=1$. For  $e\in \{1, 2\}$, the values of the normalized residues $\rho_i\equiv F_i^e\pmod{F_j}$ with $0\le \rho_i<F_j-1$ are obtained.
\medskip

\end{abstract} 

\section{Introduction}
\label{s introduction}
The Fibonacci sequence $\mathcal{F}=(F_i:i\ge 0)$ is defined by recursion:  $F_0=0$, $F_1=1$ and $F_{i+1}=F_i+F_{i-1}$ for all $i\ge 1$. Let $m\ge 3$ be an integer. The periodicity of the sequence ${\mathcal F}(m)=(F_i \pmod{m}: i\ge 0)$ has attracted a lot of attention, see for example M.~Renault~\cite{Renault1,Renault2} and the references cited therein. In particular, if $m=F_j$ is a Fibonacci number with $j\ge 4$, A.~Ehrlich~\cite{Ehrlich} proves that the sequence 
$\mathcal{F}(F_j)$ has periodicity either $2j$ if $j$ is even or periodicity $4j$ if $j$ is odd. 
Our aim in the present paper is to study the periodicity of the power Fibonacci sequences 
$$
{\mathcal F}^e(F_j)=(F_i^e \!\!\!\!\pmod{F_j}: i\ge 0),
$$
for integers $j\ge 0$ and $e\ge 1$.

Throughout, we fix representatives of the classes modulo $F_j$. We set $\rho_i(F_j,e)$ defined by $\rho_i(F_j,e)\equiv F_i^e\pmod{F_j}$ with $0\le \rho_i(F_j,e)\le F_j-1$. 
If $F_j$ and $e$ are clear from the context, we shall write simply $\rho_i$ instead of $\rho_i(F_j,e)$. 


The product of two sequences ${\mathcal F}^{e_1}(F_j)$ and ${\mathcal F}^{e_2}(F_j)$ is the product term to term modulus $F_j$, that is, ${\mathcal F}^{e_1}(F_j)\cdot {\mathcal F}^{e_2}(F_j)={\mathcal F}^{e_1+e_2}(F_j)$. 
A sequence $S=(a_i:i\ge 0)$ is \emph{periodic} if it exists an integer $p\ge 1$, called a \emph{period}, such that $a_{i+p}=a_i$ for all $i\ge 0$. In this case, the sequence $S$ is denoted by $S=[a_0,a_1,\ldots,a_{p-1}]$. It is clear that if $p$ is a period of a sequence $S$, then $qp$ is also a period of $S$ for all integers $q\ge 1$. If ${\mathcal F}^e(F_j)$ is periodic, the minimum period is called its \emph{periodicity} and will be denoted by $\pi(F_j,e)$. In this case, it is clear that $\pi(F_j,qe)\mid \pi(F_j,e)$ for all integers $q\ge 1$, and that $\pi(F_j,e)$ is a divisor of all periods of ${\mathcal F}^e(F_j)$.

If $j=0$, we have $F_j=0$ and the sequence ${\mathcal F}^e(F_j)$ is the sequence $(F_i^e:i\ge 0)$. Since 
$F_i^e<F_{i+1}^e$ for all $i\ge 2$, the sequence ${\mathcal F}^e(F_0)$ is not periodic.  

If $j=1$ or $j=2$, we have $F_j=1$. Then, ${\mathcal F}^e(F_1)={\mathcal F}^e(F_2)=[0]$ and $\pi(F_1,e)=\pi(F_2,e)=1$.

If $j=3$, we have $F_j=2$. Then ${\mathcal F}(F_3)=[0,1,1]$ and  ${\mathcal F}^e(F_3)=[0,1,1]$ for all $e\ge 1$. Thus, $\pi(F_3,e)=3$. 

From now on, we assume that $j\ge 4$.

In Sections~\ref{s e=1} and \ref{s e=2} we give the exact values of $\rho_i(F_j,1)$ and $\rho_i(F_j,2)$, respectively, and the periods $\pi(F_j,1)$ and $\pi(F_j,2)$. In  Section~\ref{s e=e} we obtain the periodicity of ${\mathcal F}^e(F_j)$ for any $e>2$.

Along the paper we use the following well known properties of Fibonacci numbers (the proofs can be found in~D.~M.~Burton~\cite{Burton}, T.~Koshly~\cite{Koshy}, or N.~N.~Vorob'ev~\cite{Vorobev}):
\begin{align} 
& \gcd(F_n,F_m)=F_{\gcd(n,m)} \quad \mbox{for all $m,n\ge 0$, $(m,n)\ne(0,0)$};\label{gcd}\\
& F_{n+m}=F_{n-1}F_m+F_{n}F_{m+1}\quad \mbox{for all $m\ge 0$ and $n\ge 1$}; \label{m+n}\\
& F_n^2-F_{n-r}F_{n+r}=(-1)^{n-r}F_r^2\quad \mbox{for all $n\ge r\ge 0$}\quad \mbox{(Catalan identity)};\label{catalan}\\
&F_n^2-F_{n-1}F_{n+1}=(-1)^{n-1}\quad \mbox{for all $n\ge 1$}\quad \mbox{(Cassini identity)}.\label{cassini}
\end{align}
Let $\phi=(1+\sqrt{5})/2$ and $\psi=(1-\sqrt{5})/2$. 
In the  application of Carmichael's Theorem in Section~\ref{s e=e} we also need the well known fact that for all integers $n\ge 0$, it holds
$$
F_n=(\phi^n-\psi^n)/\sqrt{5}=(\phi^n-\psi^n)/(\phi-\psi).
$$

Finally, a remark about notation. Throughout,  if the modulus is not explicitly  indicated, in all congruences it is understood that this is the Fibonacci number $F_j$ from the context.

\section{The sequences ${\mathcal F}(F_j)$}
\label{s e=1}

A.~Ehrlich~\cite{Ehrlich} proved that if $j\ge 4$ then $\pi(F_j,1)=2j$ if $j$ is even and $\pi(F_j,1)=4j$ if $j$ is odd. In this section, we calculate explicitly the terms of ${\mathcal F}(F_j)$ and obtain the same result.

The cases $j\in\{4,5,6,7\}$ are easily obtained:
\begin{equation}
\label{j4567}
\begin{array}{lr@{\hskip 2pt}l}
F_4=3,  & {\mathcal F}(F_4)=&(0,1,1,2,0,2,2,1,0,1,\ldots)\\
        &       =&[0,1,1,2,0,2,2,1]. \\
F_5=5,  & {\mathcal F}(F_5)=&(0,1,1,2,3,0,3,3,1,4,0,4,4,3,2,0,2,2,4,1,0,1,\ldots)\\
        &       =&[0,1,1,2,3,0,3,3,1,4,0,4,4,3,2,0,2,2,4,1]. \\
F_6=8,  & {\mathcal F}(F_6)=&(0,1,1,2,3,5,0,5,5,2,7,1,0,1,\ldots)\\
        &       =&[0,1,1,2,3,5,0,5,5,2,7,1].\\
F_7=13, & {\mathcal F}(F_7)=&(0,1,1,2,3,5,8,0,8,8,3,11,1,12,0,12,12,11,10,8,5,0,5,5,10,2,12,1,0,1,\ldots) \\
        &       =&[0,1,1,2,3,5,8,0,8,8,3,11,1,12,0,12,12,11,10,8,5,0,5,5,10,2,12,1].
\end{array}
\end{equation}

Then,
\begin{equation}
\label{p4567}
\pi(F_4,1)=8, \quad \pi(F_5,1)=20,\quad \pi(F_6,1)=12, \quad \pi(F_7,1)=28.
\end{equation}
 
\begin{thm} 
\label{thm e=1}
Let $j\ge 4$ be an integer. 
\begin{enumerate}[\rm (i)]
\item If $j$ is even, then $\pi(F_j,1)=2j$ and the values of $\rho_i=\rho_i(F_j,1)$ for $i\in\{0,\ldots,2j-1\}$ are
\begin{enumerate}
\item[\rm (i.1)] $\rho_i=F_i$  \  if \  $i<j$;
\item[\rm (i.2)]  $\rho_i=0$ \ if \   $i=j$;
\item[\rm (i.3)]  $\rho_i=F_{2j-i}$ \ if \  $j<i<2j$\ and\ $i$\ is odd; 
\item[\rm (i.4)]  $\rho_i=F_j-F_{2j-i}$ \ if \  $j<i<2j$\ and\  $i$\ is even.
\end{enumerate}
\item If $j$ is odd, then $\pi(F_j,1)=4j$ and the values of $\rho_i=\rho_i(F_j,1)$ for $i\in\{0,\ldots,4j-1\}$ are
\begin{enumerate}[\rm 1)]
\item[\rm (ii.1)] $\rho_i=F_i$ \ if \  $i<j$; 
\item[\rm (ii.2)]  $\rho_i=0$ \ if \ $i=j$;
\item[\rm (ii.3)]  $\rho_i=F_{2j-i}$ \  if \ $j<i<2j$ and $i$ is even;
\item[\rm (ii.4)]  $\rho_i=F_j-F_{2j-i}$ \ if \  $j<i<2j$ and $i$ is odd;
\item[\rm (ii.5)]  $\rho_i=0$ \ if \ $i=2j$; 
\item[\rm (ii.6)]  $\rho_i=F_j-F_{i-2j}$ \  if \ $2j<i<3j$; 
\item[\rm (ii.7)]  $\rho_i=0$ \ if \ $i=3j$; 
\item[\rm (ii.8)]  $\rho_i=F_j-F_{4j-i}$ \ if \ $3j<i<4j$  and $i$ is even;
\item[\rm (ii.9)]  $\rho_i=F_{4j-i}$ \ if \ $3j<i<4j$ and $i$ is odd.
\end{enumerate}
\end{enumerate}
\end{thm}
\begin{pf}
By~(\ref{j4567}) and~(\ref{p4567}), the result is true for $j\in\{4,5,6,7\}$.
Now, let $j\ge 8$ and assume the properties are true for values $<j$. First we check that the values of $\rho_i$ for $i\in\{0,\ldots,2j-1\}$ when $j$ is even and for $i\in\{0,\ldots,4j-1\}$ when $j$ is odd are those given in the respective list of items. 

(i)  Assume that $j$ is even. Then, 

(i.1) and (i.2), are obvious. 
 
(i.3) and (i.4). By induction on $i\in\{j+1,\ldots,2j-1\}$. For $i=j+1$ we have $i$ odd,  $2j-i=2j-(j+1)=j-1$, and
$$
\rho_{j+1}\equiv F_{j+1}=F_j+F_{j-1}\equiv F_{j-1}=F_{2j-i}.
$$
For $i=j+2$, we have $i$ even, $2j-i=2j-(j+2)=j-2$, and
$$
\rho_{j+2}\equiv F_{j+2}=F_{j+1}+F_j=F_j+F_{j-1}+F_j\equiv F_{j-1}=F_j-F_{j-2}=F_j-F_{2j-i}.
$$
If $i\ge j+3$ and $i$ is odd, 
$$
\rho_i\equiv F_{i}
=F_{i-1}+F_{i-2}
\equiv \rho_{i-1}+\rho_{i-2}
=F_j-F_{2j-i+1}+F_{2j-i+2}
\equiv F_{2j-i+2}-F_{2j-i+1} 
=F_{2j-i}.
$$
If $i\ge j+4$ and $i$ is even,
$$
\rho_i\equiv F_i=F_{i-1}+F_{i-2}
\equiv \rho_{i-1}+\rho_{i-2}
=F_{2j-i+1}+F_j-F_{2j-i+2}
=F_j-(F_{2j-i+2}-F_{2j-i+1})
=F_j-F_{2j-i}.
$$

For $1\le i\le 2j-1$, the unique value $i$ such that $\rho_i=0$ is $i=j$, but $\rho_{j+1}=F_{j-1}\ne 1$.
Hence, ${\mathcal F}(F_j)$ has not a period $\le 2j-1$. Now, $\rho_{2j-2}=F_j-F_2=F_j-1\equiv -1$ and $\rho_{2j-1}=F_1=1$. 
Hence, $\rho_{2j}=0$ and $\rho_{2j+1}=1$. Thus, $\pi(F_j,1)=2j$. 

(ii) Now, assume that $j$ is odd. Then,

(ii.1) and (ii.2) are obvious.

(ii.3) and (ii.4) By induction on $i$. If $i=j+1$, then $i$ is even, $2j-i=j-1$, and 
$$
\rho_{j+1}\equiv F_{j+1}=F_j+F_{j-1}\equiv F_{j-1}=F_{2j-i}.
$$
If $i=j+2$, then $i$ is odd , $2j-i=j-2$, and
$$
\rho_{j+2}\equiv F_{j+2}=F_{j+1}+F_j\equiv F_{j+1}\equiv \rho_{j+1}=F_{j-1}=F_j-F_{j-2}=F_j-F_{2j-i}.
$$
If $j+3\le i<2j-1$ and is $i$ even, then
$$
\rho_i\equiv F_i=F_{i-1}+F_{i-2}\equiv \rho_{i-1}+\rho_{i-2}=F_j-F_{2j-i+1}+F_{2j-i+2}\equiv F_{2j-i}.
$$
If $j+3\le i<2j-1$ and is $i$ is odd, then
$$
\rho_i\equiv F_i=F_{i-1}+F_{i-2}\equiv F_{2j-i+1}+F_j-F_{2j-i+2}=F_j-(F_{2j-i+2}-F_{2j-i+1})=F_j-F_{2j-i}.
$$

(ii.5) Since $\rho_{2j-2}=F_2=1$ and $\rho_{2j-1}=F_j-F_{1}=F_j-1$, we have $\rho_{2j}\equiv F_j\equiv 0$. 

(ii.6) By induction on $i$. If $i=2j+1$, we have $i-2j=1$ and 
$$
\rho_{2j+1}\equiv F_{2j+1}=F_{2j}+F_{2j-1}\equiv \rho_{2j}+\rho_{2j-1}=\rho_{2j-1}=F_j-F_{1}=F_j-F_{i-2j}.
$$
If $i=2j+2$, we have $i-2j=2$ and 
$$
\rho_{2j+2}\equiv F_{2j+2}=F_{2j+1}+F_{2j}\equiv F_j-F_1=F_j-F_2=F_j-F_{i-2j}.
$$
If $2j+2<i<3j$, then
$$
\rho_i
\equiv F_i=F_{i-1}+F_{i-2}
\equiv F_j-F_{i-1-2j}+F_j-F_{i-2-2j}
\equiv F_j-(F_{i-1-2j}+F_{i-2-2j})
=F_j-F_{i-2j}.
$$

(ii.7) Since $\rho_{3j-2}=F_j-F_{j-2}$ and $\rho_{3j-1}=F_j-F_{j-1}$, we have $\rho_{3j}\equiv -F_{j-1}-F_{j-2}=-F_j\equiv 0$.

(ii.8) and (ii.9) By induction on $i$. If $i=3j+1$, then $i$ is even, $4j-i=j-1$, and 
$$
\rho_{3j+1}\equiv F_{3j+1}= F_{3j}+F_{3j-1}\equiv F_{3j-1}\equiv \rho_{3j-1}
=F_j-F_{j-1}=F_j-F_{4j-i}.
$$
If $i=3j+2$, then $i$ is odd, $4j-i=j-2$ and
$$
\rho_{3j+2}\equiv F_{3j+1}+F_{3j}
\equiv F_j-F_{j-1}+F_{3j}
=F_j-F_{j-1}=F_{j-2}=F_{4j-i}.
$$
If $3j+2<i<4j$ and $i$ is even, then
$$
\rho_i\equiv F_i=F_{i-1}+F_{i-2}=F_{4j-i+1}+F_j-F_{4j-i+2}=F_j-(F_{4j-i+2}-F_{4j-i+1})=F_j-F_{4j-i}.
$$
If $3j+2<i<4j$ and $i$ is odd, then
$$
\rho_i\equiv F_i=F_{i-1}+F_{i-2}=F_j-F_{4j-i+1}+F_{4j-i+2}\equiv F_{4j-i+2}-F_{4j-i+1}=F_{4j-i}.
$$
 
The unique values of $i\in\{1,\ldots,4j-1\}$ such that $\rho_i=0$ are $i=j$, $i=2j$ and $i=3j$.
But $\rho_{j+1}=F_{j-1}\ne 1$, $\rho_{2j+1}=F_j-F_1\ne 1$ and $\rho_{3j+1}=F_j-F_{j-1}=F_{j-2}\ne 1$.
Therefore, ${\mathcal F}(F_j)$ has no period $<4j$. 

Now, 
$$
\rho_{4j}\equiv F_{4j}=F_{4j-1}+F_{4j-2}\equiv F_1+F_j-F_2\equiv 0,
$$
and 
$$
\rho_{4j+1}\equiv F_{4j+1}=F_{4j}+F_{4j-1}\equiv F_{4j-1}\equiv F_1=1.
$$
Hence $\pi(F_j,1)=4j$. \qed
\end{pf}

\section{The square sequences ${\mathcal F}^2(F_j)$}
\label{s e=2}

In this section we calculate the terms of ${\mathcal F}^2(F_j)$ and prove that $\pi(F_j,2)=j$ if $j$ is even and $\pi(F_j,2)=2j$ if $j$ is odd.

\begin{lem} 
\label{k+a,k-a}
Let $k\ge 2$ be an integer, and $\alpha\in \{0,1,\ldots,k\}$. It holds 
\begin{enumerate}[\rm (i)] 
\item $F_k^2<F_{2k}$;
\item $F^2_{k+\alpha}\equiv F^2_{k-\alpha} \pmod{F_{2k}}$;
\item $F_{k+1}^2<F_{2k+1}$;
\item $F^2_{k+1+\alpha}\equiv -F_{k-\alpha}^2 \pmod{F_{2k+1}}$.
\end{enumerate}
\end{lem}
\begin{pf}
(i)  Take $m=n=k$ in~(\ref{m+n}). Then,   $F_{2k}=F_{k-1}F_k+F_kF_{k+1}>F_kF_{k+1}>F_k^2$.

(ii) Take $n=k+\alpha$ and $r=k-\alpha$ in~(\ref{catalan}). Then, $n-r=2\alpha$, $n+r=2k$, and
      $$
      F_{k+\alpha}^2-F_{2\alpha}F_{2k}=(-1)^{2\alpha}F_{k-\alpha}^2=F_{k-\alpha}^2.
      $$
       By taking modulus $F_{2k}$ we get $F_{k+\alpha}^2\equiv F^2_{k-\alpha}\pmod{F_{2k}}$. 

(iii) Take $m=k$ and $n=k+1$ in~(\ref{m+n}). Then, one has
$$
F_{2k+1}=F_{k+(k+1)}
=F_{k-1}F_{k+1}+F_kF_{k+2}
>F_{k-1}F_{k+1}+F_kF_{k+1}
=(F_{k-1}+F_k)F_{k+1}=F_{k+1}^2.
$$

(iv) Take $n=k+1+\alpha$ and $r=k-\alpha$ in~(\ref{catalan}). We have
 $n-r=2\alpha+1$, $n+r=2k+1$, and
 $$
 F_{k+\alpha+1}^2-F_{2\alpha+1}F_{2k+1}=(-1)^{2\alpha+1}F_{k-\alpha}^2=-F_{k-\alpha}^2.
 $$
 Thus, $F_{k+\alpha+1}^2\equiv -F_{k-\alpha}^2 \pmod{F_{2k+1}}$.\qed
\end{pf} 

\begin{thm}   
\label{e=2}
Let $j\ge 4$ be an integer and $\rho_i=\rho_i(F_j,2)$.
\begin{enumerate}[\rm (i)] 
\item If $j=2t$ is even, then ${\mathcal F}^2(F_j)$ has periodicity $j$ and the values of $\rho_i$ for $i\in\{0,\ldots,j-1\}$ are
\begin{enumerate}
\item[\rm (i.1)] \label{qeven}
$\rho_i=F_i^2$ \ if \  $i\in\{0,\ldots,t\}$;
\item[\rm (i.2)] \label{odd}
$\rho_i=F_{j-i}^2$ \ if \ $i\in\{t+1,\ldots,j-1\}$.
\end{enumerate}
\item If $j=2t+1$ is odd, then ${\mathcal F}^2(j)$ has periodicity $2j$ and the values of $\rho_i$ for $i\in\{0,\ldots,2j-1\}$ are
\begin{enumerate}
\item[\rm (ii.1)]
$\rho_i=F_i^2$ \  if \  $i\in\{0,\ldots,t+1\}$;
\item[\rm (ii.2)]
$\rho_i=F_j-F_{j-i}^2$ \  if \ $i\in\{t+2,\ldots,j-1\}$;
\item[\rm (ii.3)]
$\rho_i=0$ \ if \ $i=j$;
\item[\rm (ii.4)]
$\rho_i=\rho_{2j-i}$ \ if \  if  $i\in\{j+1,\ldots,2j-1\}$. 
\end{enumerate}
\end{enumerate}
\end{thm}
\begin{pf}
The proof is very similar to that of Theorem \ref{thm e=1}.
As for (i.1), if $i\in\{0,\ldots,t\}$, by Lemma~\ref{k+a,k-a}(i), we have $\rho_i\equiv F_i^2<F_{2i}\le F_{2t}=F_j$. Hence, $\rho_i=F_i^2$. 

(i.2) For $i\in\{t+1,\ldots,j-1\}$, set $\alpha=i-t$. Then, $j-i=2t-(t+\alpha)=t-\alpha$. By Lemma~\ref{k+a,k-a}(ii), we have $F_{t+\alpha}^2\equiv F_{t-\alpha}^2 \pmod{F_{2t}}$. Thus, 
$$
\rho_i\equiv F_i^2=F_{t+\alpha}^2\equiv F_{t-\alpha}^2=F_{j-i}^2.
$$
Since $j-i\le t$, we have $F^2_{j-i}<F_j$. Hence, $\rho_i=F_{j-i}^2$.

Now, we prove that $\pi(F_j,2)=j$. Indeed,
from~(\ref{cassini}), we have $F_{j+1}^2-F_jF_{j+2}=(-1)^{j}=1$. Hence, $F_{j+1}^2\equiv 1$. Then, by (\ref{m+n}), 
$$
F_{i+j}^2=(F_{i-1}F_j+F_iF_{j+1})^2\equiv F_i^2F_{j+1}^2\equiv F_i^2.
$$ 
This implies that $\pi(F_j,2)\mid j$. But, for $i\in\{1,\ldots,j-1\}$, we have $\rho_i\ne 0$. 
Hence, $\pi(F_j,2)=j$. 

(ii.1) For $i\in\{0,1,\ldots,t+1\}$, by Lemma~\ref{k+a,k-a}(iii), we have $0\le F_i^2\le F_{t+1}^2<F_{2t+1}=F_j$. Hence, $\rho_i=F_i^2$.

(ii.2) For $i\in\{t+2,\ldots,j-1\}$, define $\alpha$ by $i=t+1+\alpha$. Then, $j-i=2t+1-(t+1+\alpha)=t-\alpha$ and, by  Lemma~\ref{k+a,k-a}(iv),  we have 
$$
F_i^2=F_{t+\alpha+1}^2\equiv -F_{t-\alpha}^2=-F_{j-i}^2\equiv F_j-F_{j-i}^2.
$$
Since $F_{j-i}^2<F_j$ due to $j-i\le t+1$, we conclude
$\rho_i=F_j-F_{j-i}^2$.

(ii.3) Obviously, $\rho_j=0$. 

(ii.4) For $i\in\{j+1,\ldots,2j\}$, define $\alpha$ by $i=j+\alpha$. Again, by Catalan's identity~(\ref{catalan}) 
with $r=\alpha$, we have
$$
F_{j+\alpha}^2-F_jF_{j+2\alpha}=(-1)^jF_{\alpha}^2=-F_\alpha^2
$$
so that $F_{j+\alpha}^2\equiv -F_{\alpha}^2$. Thus, for $\alpha\in\{1,\ldots,t+1\}$, one has
$$
\rho_i=\rho_{j+\alpha}\equiv -F_{\alpha}^2\equiv F_j-F_{\alpha}^2=\rho_{j-\alpha}=\rho_{2j-i}.
$$

Finally, we prove that $\pi(F_j,2)=2j$. 
By~(\ref{cassini}) we have 
$$
F_{2j+1}^2-F_{2j}F_{2j+2}=(-1)^{2j}=1.
$$
Since $F_{2j}\equiv 0$, we have $F_{2j+1}^2\equiv 1$. 
By~(\ref{m+n}), it holds
$$
F_{i+2j}^2=(F_{i-1}F_{2j}+F_iF_{2j+1})^2 \equiv F_i^2F_{2j+1}^2\equiv F_i^2.
$$
Thus, the periodicity of ${\mathcal F}^2(F_j)$ is at most, $2j$.
Now, for $i\in\{1,\ldots,2j-1\}$, the unique $i$ with $\rho_i=0$ is $i=j$. Since $\rho_{j+1}=F_j-F_1^2=F_j-1\ne 1$,
we have $\pi(F_j,2)>j$. We conclude $\pi(F_j,2)=2j$. \qed
\end{pf}

As examples of this, take $j=6$ so that $F_6=8$, $t=4$, and
$$
\begin{array}{ll}
 {\mathcal F}(F_6)  &=[0,1,1,2,3,5,0,5,5,2,7,1], \\[4pt]
 {\mathcal F}^2(F_6)&=[0,1,1,4,1,1].
\end{array}
$$
For $j=7$ we have $F_7=13$, $t=6$, and 
$$
\begin{array}{ll}
{\mathcal F}(F_7)
&=[0,1,1,2,3,\phantom{1}5,\phantom{1}8,0,\phantom{1}8,\phantom{1}8,3,11,1,12,0,12,12,11,10,8,5,0,5,5,10,2,12,1], \\[3pt]
{\mathcal F}^2(F_7)
&=[0,1,1,4,9,12,12,0,12,12,9,\phantom{1}4,1,\phantom{1}1].
\end{array}
$$

\section{The power sequences ${\mathcal F}^e(F_j)$}  
\label{s e=e}

Consider first the special case $j=6$.

\begin{prop}
\label{j=6}
\begin{enumerate}[\rm (i)]
\item If $e\ge 1$ is odd, then $\pi(F_6,e)=12$;
\item $\pi(F_6,2)=6$, and if $e\ge 4$ is even, then $\pi(F_6,e)=3$.
\end{enumerate}
\end{prop}
\begin{pf}
(i) Since $F_6=8$, a direct calculation gives
\begin{align*}
{\mathcal F}(F_6)&=[0, 1, 1, 2, 3, 5, 0, 5, 5, 2, 7, 1],\\
{\mathcal F}^2(F_6)&=[0, 1, 1, 4, 1, 1],\\
{\mathcal F}^3(F_6)&=[0, 1, 1, 0, 3, 5, 0, 5, 5, 0, 7, 1],\\
{\mathcal F}^4(F_6)&=[0, 1, 1, 0, 1, 1, 0, 1, 1, 0, 1, 1]\\
        &=[0,1,1].
\end{align*}
It is clear that ${\mathcal F}^{4\alpha}(F_6)=[0,1,1]$ for all integers $\alpha\ge 1$. Note that $F_i^4\equiv 0$ if and only if $i\equiv 0\pmod{3}$, and if $F_i\not\equiv F_i^3$ then $i\equiv 0\pmod{3}$. Thus,
$$
{\mathcal F}^{4\alpha+1}(F_6)={\mathcal F}^{4\alpha}(F_6){\mathcal F}(F_6)={\mathcal F}^{4\alpha}(F_6){\mathcal F}^3(F_6)={\mathcal F}^{4\alpha+3}(F_6).
$$
If $e$ is odd, then $e=4\alpha+1$ or $e=4\alpha+3$ for some $\alpha\ge 1$ and we have
\begin{align*}
{\mathcal F}^{4\alpha+1}(F_6)={\mathcal F}^{4\alpha+3}(F_6)
&=[0,1,1]\cdot [0, 1, 1, 0, 3, 5, 0, 5, 5, 0, 7, 1]\\
&=[0,1,1,0,3,5,0,5,5,0,7,1].
\end{align*}
Hence, for $e$ odd, we get $\pi(F_6,e)=12$.

(ii) It is clear that $\pi(F_6,2)=6$. If $e\ge 4$ is even, then $e=4\alpha$ or $e=4\alpha+2$ for some $\alpha\ge 1$. 
 Since ${\mathcal F}^{4\alpha}(F_6)=[0,1,1]$, we have
 $$
 {\mathcal F}^{4\alpha+2}(F_6)={\mathcal F}^{4\alpha}(F_6)\cdot {\mathcal F}^2(F_6)=[0,1,1]\cdot[0, 1, 1, 4, 1, 1]=[0,1,1].
 $$
 Hence $\pi(F_6,e)=3$ for all $e\ge 4$ even.\qed
\end{pf}

Now we discuss the positions of $0$ in the sequence ${\mathcal F}^e(F_j)$. To this end, we use the Carmicheal Theorem (see \cite{Carmichael} and \cite{Carmichael2}, or \cite{Yabuta}) applied to Fibonacci numbers:

\begin{thm}[Carmichael] 
\label{Carmichael} 
Let $j$ be an integer with $j\ge 3$ and $j\ne 12$. Then, there exists a prime $p$ such that $p\mid F_j$ but $p\nmid F_i$ for all $i\in\{1,\ldots,j-1\}$.
\end{thm}

\begin{prop}
\label{pos0, e>2}\
Let $j$ and $e$ be integers such that $e\ge 1$ and $4\le j\ne 6$. Then $F_i^e\equiv 0\pmod{F_j}$ if and only if $i\equiv 0\pmod{j}$.
\end{prop} 

\begin{pf}
If $i\equiv 0\pmod{j}$, then $\gcd(i,j)=j$, which implies $\gcd(F_i,F_j)=F_j$.
Hence $F_i\equiv 0$ and $F_i^e\equiv 0$ for all $e\ge 1$.

Conversely, assume $i\not\equiv 0\pmod{j}$.
Consider first the special case $F_{12}=144=2^4\cdot 3^2$. The prime factorizations of $F_i$ for $i\in\{1,\ldots,11\}$ are the following:
$$
\begin{array}{r|lllllllllll}
i & 1 & 2 & 3 & 4 & 5 & 6     & 7  & 8           & 9            & 10           & 11 \\
\hline
F_i & 1 & 1 & 2  & 3 & 5 & 8=2^3 & 13 & 21=3\cdot 7 & 34=2\cdot 17 & 55=5\cdot 11 & 89
\end{array}
$$
No $F_i$ in the table has both factors $2$ and $3$. Thus, $F_i^e$ has not both factors $2$ and $3$. Hence $F_i^e\not\equiv 0\pmod{F_{12}}$.

Suppose $j\ne 12$. By Carmichael's Theorem there exists a prime factor $p$ of $F_j$ which is not a prime factor of  $F_i$ for $i\in\{1,\ldots,j-1\}$. Then $p$ is not a prime factor of $F_i^e$, and $F_i^e\not\equiv 0$. 

If $i=qj+r$ with $q\ge 1$ and $1\le r\le j-1$, we have $\gcd(F_j,F_{qj})=F_j$. Then,
$F_{qj}\equiv 0$ and, by (\ref{m+n}),
$$
F_i=F_{qj+r}=F_{qj-1}F_r+F_{qj}F_{r+1}\equiv F_{qj-1}F_r.
$$
Since $\gcd(qj-1,j)=1$, it follows that $F_{qj-1}$ is a unit modulo $F_j$ and thus $F_{qj-1}^e$ is a unit too. Since $1\le r<j-1$, by the previous paragraph we have $F_r^e \not\equiv 0$. Thus,
$$
F_i^e\equiv F_{qj-1}^eF_r^e\not\equiv 0. \qed
$$
\end{pf}

\begin{thm}
\label{main}
Let $j$ and $e$ be integers such that $j\ge 4$, $j\ne 6$, and $e\ge 1$. 
\begin{enumerate}[\rm (i)]
\item Assume that $j$ is even.
\begin{enumerate}
\item[\rm (i.1)] If $e$ is even, then $\pi(F_j,e)=j$;
\item[\rm (i.2)] if $e$ is odd, then $\pi(F_j,e)=2j$.
\end{enumerate}
\item Assume that $j$ is odd.
\begin{enumerate}
\item[\rm (ii.1)] If $e\equiv 0\pmod{4}$, then $\pi(F_j,e)=j$.
\item[\rm (ii.2)] if $e\equiv 2\pmod{4}$, then $\pi(F_j,e)=2j$.
\item[\rm (ii.3)] if $e$ is odd, then $\pi(F_j,e)=4j$.
\end{enumerate}
\end{enumerate}
\end{thm}
\begin{pf}
Due to Proposition~\ref{pos0, e>2}, we know that $\pi(F_j,e)$ is a multiple of $j$. 

(i.1) Define $\alpha$ by $e=2\alpha$. 
We shall see that $F_{i+j}^{2\alpha}\equiv F_i^{2\alpha}$ by induction on $\alpha$.
For $\alpha=1$, the result is true by Theorem~\ref{e=2}. If $\alpha\ge 2$ and it is true for $\alpha-1$, then
$$
F_{i+j}^e\equiv F_{i+j}^{2\alpha}
=F_{i+j}^{2(\alpha-1)}F_{i+j}^2
\equiv F_i^{2(\alpha-1)}F_i^2=F_i^{2\alpha}=F_i^e.
$$

(i.2) Define $\alpha$ by $e=2\alpha+1$. 
To see that  $\pi(F_j,e)\ge 2j$ it is sufficient to show that $F_{j+1}^e\not\equiv 1$. By Theorem~\ref{e=2}, we have 
$$
F_{j-1}^{2\alpha}=(F_{j-1}^2)^\alpha=F_1^\alpha=1,
$$
and
$$
F_{j+1}^e
=F_{j+1}^{2\alpha+1}
=(F_j+F_{j-1})^{2\alpha+1}
\equiv F_{j-1}^{2\alpha+1}
\equiv F_{j-1}^{2\alpha}F_{j-1}
\equiv F_{j-1}\not\equiv 1.
$$
Thus, $\pi(F_j,e)\ge 2j$. We shall see that $F_{i+2j}^{2\alpha+1}\equiv F_i^{2\alpha+1}$ by induction on $\alpha$. For $\alpha=0$ the result is true by Theorem~\ref{thm e=1}. Assume  $\alpha\ge 1$ and the result true for $\alpha-1$. By (i.1) we have $F_{i+2j}^{2\alpha}\equiv F_i^{2\alpha}$. By Theorem~\ref{thm e=1}, we have $F_{i+2j}\equiv F_i$. Then,
$$
F_{i+2j}^{2\alpha+1}=F_{i+2j}^{2\alpha}F_{i+2j}\equiv F_i^{2\alpha}F_i=F_i^{2\alpha+1}.
$$ 

(ii) Let $j=2t+1$. 

(ii.1) Let $e=4\alpha$. It is enough to prove that $j$ is a period of ${\mathcal F}^e(F_j)$. 

Consider first the case $\alpha=1$; that is, $e=4$. By the identity~(\ref{m+n}), we have
$$
0\equiv F_j=F_{2t+1}=F_{(t+1)+t}=F_t^2+F_{t+1}^2.
$$
Thus, $F_t^4\equiv F_{t+1}^4$. Then, $\rho_t(F_j,4)=\rho_{t+1}(F_j,4)$. By Theorem~\ref{e=2}, for 
$i\in\{t+2,\ldots,j-1\}$, we have
$F_i^4\equiv (F_j-F_{j-i}^2)^2\equiv F_{j-i}^4$. Thus, the first $j$ terms of the sequence ${\mathcal F}^4(F_j)$ are
$$
F_0^4,F_1^4,\ldots,F_{t-1}^4,F_t^4,F_t^4,F_{t-1}^4,\ldots,F_1^4.
$$
Since for $i\in\{j+1,$\ldots$,2j-1\}$ we have $\rho_i(F_j,4)\equiv \rho_i(F_j,2)^2=\rho_{2j-i}(F_j,2)^2=\rho_{2j-i}(F_j,4)$, we see that $j$ is a period of $\mathcal{F}^4(F_j)$. Hence $\pi(F_j,4)=j$. This implies $\pi(F_j,4\alpha)=j$.

(ii.2) Let $e=4\alpha+2$. It is enough to prove that $j$ is not a period. We shall see that
$F_{j+1}^e\not\equiv F_1$.  By Theorem~\ref{e=2}, one has 
$$ 
F_{j+1}^2\equiv \rho_{j-1}(F_j,2)=F_j-F_1^2=F_j-1\equiv -1.
$$
By (ii.1), it follows $F_{j+1}^{4\alpha}\equiv F_1^{4\alpha}\equiv 1$. Then, 
$$
F_{j+1}^{4\alpha+2}=F_{j+1}^{4\alpha}F_{j+1}^2\equiv -1\not\equiv F_1.
$$
Hence $j$ is not a period of ${\mathcal F}^e(F_j)$. Therefore, we get $\pi(F_j,e)=2j$.

(ii.3) Let $e=2\alpha+1$. Since $\pi(F_j,1)=4j$, we have $\pi(F_j,e)\in\{j,2j,3j,4j\}$. 
By (ii.1), it follows that $F_{j+1}^{4\alpha}\equiv F_1^{4\alpha}=1$ and, by Theorem~\ref{thm e=1}, we get 
$F_{j+1}=F_{j-1}\not\equiv 1$. Then,
$$
F_{j+1}^{4\alpha+1}=F_{j+1}^{4\alpha}F_{j+1}\equiv F_{j+1}=F_{2j-(j+1)}=F_{j-1}\not\equiv 1\,,
$$
and hence $\pi(F_j,e)\ne j$.

Analogously, we find
$$
F_{2j+1}^{4\alpha+1}=F_{2j+1}^{4\alpha}F_{2j+1}\equiv F_{2j+1}\equiv F_j-F_{2j+1-2j}=F_j-1\not\equiv 1.
$$
so that $\pi(F_j,e)\ne 2j$.

Finally, from
$$
F_{3j+1}^{4\alpha+1}
=F_{3j+1}^{4\alpha}F_{3j+1}
\equiv F_{3j+1}
\equiv F_{4j-(3j+1)}
=F_{j-1}\not\equiv 1\,,
$$
we can conclude that $\pi(F_j,e)=4j$. \qed

\end{pf}

\end{document}